\documentclass[12pt]{article}
\usepackage[top=1in,bottom=1in,left=1in,right=1in]{geometry}
\usepackage{indentfirst}
\usepackage{amsfonts,amsmath,amsthm,amssymb,hyperref,enumerate,bm}
\usepackage{longtable}
\usepackage[caption=false]{subfig}
\usepackage{leftidx}
\usepackage{lineno}
\usepackage{color,xcolor}
\usepackage{array}
\usepackage{latexsym}
\usepackage{float}
\usepackage{color}
\usepackage{longtable,supertabular,booktabs}
\usepackage{rotating}
\usepackage{multicol,multirow}
\usepackage{epsfig}
\usepackage{fancyhdr}       
\usepackage{dsfont}
\usepackage{hyperref}
\hypersetup{
	colorlinks=true,
	linkcolor=red,
	filecolor=blue,
	urlcolor=red,
	citecolor=blue,
}

\newtheorem{theorem}{Theorem}[section]

\newtheorem{lemma}{Lemma}[section]
\newtheorem{proposition}{Proposition}[section]
\newtheorem{corollary}{Corollary}[section]
\newtheorem{remark}{Remark}[section]
\newtheorem{definition}{Definition}[section]

\newtheorem{example}{Example}[section]

\def\G1{G^\mathcal{C}}

 \newenvironment{prof}{\trivlist
      \item[\hskip\labelsep
      {\itshape Proof.}]\normalfont}
      {\hspace*{\fill}$\Box$\endtrivlist}
\begin{document}
\title{Fermat's and Catalan's equations over $M_2(\mathbb{Z})$}
\author{
Hongjian Li
\footnote{E-mail\,$:$ hongjian\_li@m.scnu.edu.cn}\qquad
Pingzhi Yuan
\footnote{Corresponding author. E-mail\,$:$ yuanpz@scnu.edu.cn.
Supported by National Natural Science Foundation of China (Grant No. 12171163).}\\
{\small\it  School of Mathematical Sciences, South China Normal University,}\\
{\small\it Guangzhou 510631, Guangdong, P. R. China} \\
}
\date{}
 \maketitle
\date{}

\noindent{\bf Abstract}\quad
Let $A=\begin{pmatrix} a & b \\ c & d \end{pmatrix}\in M_2\left(\mathbb{Z}\right)$ be a given matrix such that $bc\neq0$ and let $C(A)=\{B\in M_2(\mathbb{Z}): AB=BA\}$. In this paper, we give a necessary and sufficient condition for the solvability of the matrix equation $uX^i+vY^j=wZ^k,\, i,\, j,\, k\in\mathbb{N},\, X, \,Y,\, Z\in C(A)$, where $u,\, v,\, w$ are given nonzero integers such that $\gcd\left(u,\, v,\, w\right)=1$. From this, we get a necessary and sufficient condition for the solvability of the Fermat's matrix equation in $C(A)$. Moreover, we show that the solvability of the Catalan's matrix equation in $M_2\left(\mathbb{Z}\right)$ can be reduced to the solvability of the Catalan's matrix equation in $C(A)$, and finally to the solvability of the Catalan's equation in quadratic fields.


\medskip
\noindent{\bf MR(2020) Subject Classification} 15A24, 11D41, 11D61

\section{Introduction}
It is well-known that the Fermat's equation
$$
x^n+y^n=z^n,\quad n\in\mathbb{N}, \,n\geq3
$$
has no solutions  in positive integers $x,\, y$ and $z$ (see Wiles \cite{A}). In contrast to the classical Fermat's last theorem in integers, many scholars have studied the Fermat's equation in matrices (\cite{ch, A., A.G, A.K, M, H, L}). For example, the Fermat's equation has been investigated in rational matrices \cite{A.G}, some classes of $2\times 2$ matrices \cite{ch}, general linear group $GL_2(\mathbb{Z})$ of integral $2\times 2$ matrices with $\det=\pm1$ \cite{L} and special linear group $SL_2(\mathbb{Z})$ of integral $2\times 2$ matrices with $\det=1$ \cite{A.K}.

Another classical diophantine equation in number theory is the Catalan's equation
$$
x^m - y^n = 1,\quad m, \,n\in\mathbb{N},\, m,\, n\geq2.
$$
In 1844, Catalan \cite{cata} conjectured that this equation has no solutions in positive integers $x$ and $y$, other than the trivial solution $3^2-2^3=1$. In 2004, Mih\u{a}ilescu \cite{p.} confirmed Catalan's conjecture. In analogy with the Fermat's matrix equation, it is natural to ask whether the Catalan's equation is solvable in the ring $M_r\left(\mathbb{Z}\right)$ of all integral $r\times r$ matrices. In this paper, we will study the solvability of the Catalan's equation in $M_2(\mathbb{Z})$, i.e., $r=2$.

Let $A=\begin{pmatrix} a & b \\ c & d \end{pmatrix}\in M_2\left(\mathbb{Z}\right)$ be a given matrix such that $bc\neq0$ and let
$$C(A)=\left\{B\in M_2(\mathbb{Z}): AB=BA\right\}.$$
In this paper, we will study the solvability of the matrix equation
\begin{equation}\label{e1}
uX^i+vY^j=wZ^k,\quad i,\, j,\, k\in\mathbb{N}
\end{equation}
in $C(A)$, where $u,\, v,\, w$ are given nonzero integers such that $\gcd(u,\, v,\,w)=1$. Let $u=v=w=1$ and $i=j=k=n,\, n\geq3$. Then equation \eqref{e1} becomes the Fermat's matrix equation
$$
X^n+Y^n=Z^n,\quad n\in\mathbb{N},\, n\geq3.
$$
Let $u=w=1,\, v=-1,\, i=m,\, j=n,\, m,\, n\geq2$ and $Z=I$. Then equation \eqref{e1} becomes the Catalan's matrix equation
$$
X^m-Y^n=I, \quad m,\, n \in\mathbb{N},\, m, \,n\geq2.
$$
The rest of this paper is organized as follows. In Section \ref{s3}, we present some  properties of $C(A)$. In Section \ref{s1}, we obtain a necessary and sufficient condition for the solvability of the matrix equation \eqref{e1} in $C(A)$, and we also get a necessary and sufficient condition for the solvability of the Fermat's matrix equation in $C(A)$. In Section \ref{s2}, we study the solvability of the Catalan's matrix equation in $M_2\left(\mathbb{Z}\right)$. We show that the solvability of the Catalan's matrix equation in $M_2\left(\mathbb{Z}\right)$ can be reduced to the solvability of the Catalan's matrix equation in $C(A)$, and finally to the solvability of the Catalan's equation in quadratic fields.
\section{The properties of $C(A)$}\label{s3}
\begin{lemma}\label{l6}
Let $A=\begin{pmatrix} a & b \\ c & d \end{pmatrix}\in M_2\left(\mathbb{Z}\right)$ be a given matrix such that $bc\neq0$. Then there exists a matrix $B=\begin{pmatrix} a_1 & b_1 \\ c_1 & 0 \end{pmatrix}\in M_2\left(\mathbb{Z}\right)$, where $b_1c_1\neq0$ and $\gcd(a_1, \,b_1,\,c_1)=1$ such that $C(A)=C(B).$
\end{lemma}
\begin{prof}
Let $g=\gcd(a-d,\, b,\, c)$ and $B=\frac{1}{g}(A-dI)$. Then $B=\begin{pmatrix} (a-d)/g & b/g \\ c/g & 0 \end{pmatrix}$ is a matrix in $M_2(\mathbb{Z})$ such that
\[
\left(b/g\right)\cdot\left(c/g\right)\neq0\,\mbox{~and~}\,\gcd\left(\left(a-d\right)/g,\, b/g,\, c/g\right)=1.
\]
Note that for a matrix $C\in M_2\left(\mathbb{Z}\right)$, $AC=CA$ if and only if $BC=CB$. So $C(A)=C(B)$.
\end{prof}

By Lemma \ref{l6}, in order to study the matrix class $C(A)$, we can assume that $A=\begin{pmatrix} a & b \\ c & 0 \end{pmatrix}\in M_2(\mathbb{Z})$ is a given matrix such that $bc\neq0$ and $\gcd(a,\, b,\, c)=1$.
\begin{lemma}\label{l5}
Let $A=\begin{pmatrix} a & b \\ c & 0 \end{pmatrix}\in M_2\left(\mathbb{Z}\right)$ be a given matrix such that $bc\neq0$ and $\gcd(a,\, b, \,c)=1$. Then
$$C(A)=\left\{xI+tA: x,\, t\in \mathbb{Z}\right\}.$$
\end{lemma}
\begin{prof}
Let $B=\begin{pmatrix} a_1 & b_1 \\ c_1 & d_1 \end{pmatrix}\in C(A)$. Then $B\in M_2\left(\mathbb{Z}\right)$ and $AB=BA$. From $AB=BA$, we obtain
$$
\begin{cases}
bc_1=b_1c, \\
b(a_1-d_1)=ab_1, \\
c(a_1-d_1)=ac_1, \\
\end{cases}
$$
which imply that $c_1/c=b_1/b$ and $a_1=d_1+\frac{b_1}{b}a$. Let $c_1/c=b_1/b=p/q$,  where $p\in\mathbb{Z}$, $q\in\mathbb{N}$ and $\gcd(p,\, q)=1$. Then
\begin{equation}\label{e4}
\begin{cases}
a_1=d_1+\frac{p}{q}a, \\
b_1=\frac{b_1}{b}b=\frac{p}{q}b, \\
c_1=\frac{c_1}{c}c=\frac{p}{q}c. \\
\end{cases}
\end{equation}
From $B\in M_2\left(\mathbb{Z}\right)$ and \eqref{e4}, it follows that $q\mid\gcd(a,\, b, \,c)$, which implies that $q=1$. Therefore,
$$
B=\begin{pmatrix} d_1+pa & pb \\ pc & d_1 \end{pmatrix}=d_1I+pA.
$$
This means that $B\in\left\{xI+tA: x, \,t\in \mathbb{Z}\right\}.$

Conversely, let $B\in\left\{xI+tA: x, \,t\in \mathbb{Z}\right\}$. Then $B=xI+tA$ for some $x, \,t\in\mathbb{Z}$. Evidently, we have $B\in M_2\left(\mathbb{Z}\right)$ and $AB=BA$, so $B\in C(A)$.
\end{prof}
\begin{proposition}\label{p6}
Let $A=\begin{pmatrix} a & b \\ c & 0 \end{pmatrix}\in M_2\left(\mathbb{Z}\right)$ be a given matrix such that $bc\neq0$ and $\gcd(a,\, b, \, c)=1$. Then $C(A)$ forms a commutative ring with identity under the operations of matrix addition and multiplication.
\end{proposition}
\begin{prof}
 It is evident that $C(A)$ is a subring of $M_2(\mathbb{Z})$, and $I$ is the identity of $C(A)$. So it is sufficient to show that multiplication is commutative. Let $B_1,\, B_2\in C(A)$. Then $B_1=x_1I+t_1A$ and $B_2=x_2I+t_2A$ for some $x_1,\, x_2,\, t_1,\, t_2\in\mathbb{Z}$. Since
 $$
 \begin{aligned}
 B_1B_2&=\left(x_1I+t_1A\right)\left(x_2I+t_2A\right)=x_1x_2I+\left(x_1t_2+t_1x_2\right)A+t_1t_2A^2\\
 &=x_2x_1I+\left(x_2t_1+t_2x_1\right)A+t_2t_1A^2=\left(x_2I+t_2A\right)\left(x_1I+t_1A\right)=B_2B_1,
 \end{aligned}
 $$
it follows that  multiplication is commutative.
\end{prof}

\begin{proposition}\label{p1}
Let $A=\begin{pmatrix} a & b \\ c & 0 \end{pmatrix}\in M_2\left(\mathbb{Z}\right)$ be a given matrix such that $bc\neq0$  and $\gcd(a,\, b,\, c)=1$. Then $C(A)$ has no zero divisors if and only if $a^2+4bc$ is not a square.
\end{proposition}
\begin{prof}
Let $A_1\in C(A)$. Then $A_1=xI+tA$ for some $x,\, t\in\mathbb{Z}$. Note that $A_1=O$ if and only if $x=t=0$. We now prove necessity. Suppose that $a^2+4bc$ is a square, i.e.,
\begin{equation}\label{e2}
a^2+4bc=k^2
\end{equation}
for some integer $k$. Let $x_1=\left(-a+k\right)/2$ and $x_2=\left(-a-k\right)/2$. From \eqref{e2}, it follows that $a$ and $k$ have the same parity, which implies that $x_1$ and $x_2$ are integers. Let $B_1=x_1I+A,\, B_2=x_2I+A$. Then $B_1, B_2\in C(A)$ and $B_1\neq O, B_2\neq O$. Since
$$
\begin{aligned}
B_1B_2&=\left(x_1I+A\right)\left(x_2I+A\right)=x_1x_2I+\left(x_1+x_2\right)A+A^2\\
&=\left(x_1x_2+bc\right)I+\left(x_1+x_2+a\right)A=O,
\end{aligned}
$$
it follows that $C(A)$ has zero divisors $B_1$ and $B_2$, a contradiction.

We next prove sufficiency. Let $B$ be a nonzero element of $C(A)$. Then $B=xI+tA$ for some $x,\, t\in\mathbb{Z}$, and $x,\,t$ are not all equal to zero. Let $\lambda_1$ and $\lambda_2$ be the eigenvalues of $B$. Then
$$\lambda_{1,2}=\frac{2x+ta\pm t\sqrt{a^2+4bc}}{2}.$$
Since $a^2+4bc$ is not a square, we have $\lambda_{1}\neq0$ and $\lambda_{2}\neq0$. So $\det\left(B\right)=\lambda_1\cdot\lambda_2\neq0$. Let $B_1,\, B_2\in C(A)$ and $B_1\neq O,\, B_2\neq O$. Then $\det(B_1)\neq0$ and $\det(B_2)\neq0$. Therefore,
$$
\det\left(B_1B_2\right)=\det(B_1)\cdot\det(B_2)\neq0,
$$
which implies that $B_1B_2\neq O$. Hence, $C(A)$ has no zero divisors.
\end{prof}

\begin{corollary}\label{c0}
Let $A=\begin{pmatrix} a & b \\ c & 0 \end{pmatrix}\in M_2\left(\mathbb{Z}\right)$ be a given matrix such that $bc\neq0$ and $\gcd(a,\, b,\, c)=1$. Then $C(A)$ forms an integral domain under the operations of matrix addition and multiplication if and only if $a^2+4bc$ is not a square.
\end{corollary}
\begin{prof}
Directly from Propositions \ref{p6} and \ref{p1}.
\end{prof}

\begin{proposition}\label{p5}
Let $A=\begin{pmatrix} a & b \\ c & 0 \end{pmatrix}\in M_2\left(\mathbb{Z}\right)$ be a given matrix such that $bc\neq0$ and $\gcd(a,\, b,\, c)=1$. Then the  eigenvalues of any matrix in $C(A)$  are algebraic integers in $\mathbb{Q}\left(\sqrt{a^2+4bc}\right)$.
\end{proposition}
\begin{prof}
For any matrix $B\in C(A)$, by Lemma \ref{l5}, we have $B=\begin{pmatrix} x_0+t_0a & t_0b \\ t_0c & x_0 \end{pmatrix}$ for some $x_0,\, t_0\in\mathbb{Z}$. Then the characteristic polynomial of $B$ is
\begin{equation}\label{e5}
f(x)=x^2-\left(2x_0+t_0a\right)x+x_0^2+x_0t_0a-t_0^2bc,
\end{equation}
which is a monic polynomial with integer coefficients. From \eqref{e5}, it follows that the  eigenvalues of $B$ are
\begin{equation}\label{e6}
\frac{2x_0+t_0a\pm t_0\sqrt{a^2+4bc}}{2}\in\mathbb{Q}\left(\sqrt{a^2+4bc}\right).
\end{equation}
From \eqref{e5} and \eqref{e6}, it follows that  the  eigenvalues of $B$ are algebraic integers in $\mathbb{Q}\left(\sqrt{a^2+4bc}\right)$.
\end{prof}

\section{Matrix equation $uX^i+vY^j=wZ^k$ over $C(A)$}\label{s1}
Let $A=\begin{pmatrix} a & b \\ c & 0 \end{pmatrix}\in M_2\left(\mathbb{Z}\right)$ be a given matrix such that $bc\neq0$ and $\gcd(a,\, b,\, c)=1$, and let $u,\, v,\, w$ be given nonzero integers such that $\gcd(u,\, v,\, w)=1$. In this section, we establish a connection between the solvability of the matrix equation
\begin{equation}\label{e11}
uX^i+vY^j=wZ^k,\quad i,\, j,\, k\in\mathbb{N}
\end{equation}
in  $C(A)$  and  the solvability of the equation
\begin{equation}\label{e12}
ux^i+vy^j=wz^k,\quad i,\, j,\, k\in\mathbb{N}
\end{equation}
in  quadratic fields.

In this paper, we mainly consider the non-trivial solutions of equations \eqref{e11} and \eqref{e12}, i.e., $\det\left(XYZ\right)\neq0$ and $xyz\neq0$, respectively. Indeed, for $t=1,\, 2,\, 3$, let $X_t\in M_2\left(\mathbb{Z}\right)$ such that $\det\left(X_t\right)=tr\left(X_t\right)=0$. Then $X_t^2=O$ for $t=1,\, 2,\, 3$.  Evidently, $\left(X_1,\, X_2,\, X_3\right)$  is a solution of equation \eqref{e11} for $i,\, j,\, k\geq2$. However, these solutions are trivial.
\begin{lemma}{\rm (\cite{Xi})}\label{le5}
For a positive integer $n$, let $X=\begin{pmatrix} a & b \\ c & d \end{pmatrix}\in M_2\left(\mathbb{C}\right)$ and $
X^n=\begin{pmatrix} a_{n} & b_{n} \\ c_{n} & d_{n}\end{pmatrix}$. Assume that $x_1$ and $x_2$ are the eigenvalues of $X$. Then the following statements hold.

\begin{enumerate}
\item[\rm1)]\,  If $x_1=x_2\neq0$, then
$\left\{
  \begin{array}{l }
  a_n=\left(1+\frac{n(a-x_1)}{x_1}\right)x_1^n,\vspace{0.5ex}\\
  b_n=bnx_1^{n-1},\vspace{0.5ex} \\
  c_n=cnx_1^{n-1},\vspace{0.5ex} \\
  d_n=\left(1+\frac{n(d-x_1)}{x_1}\right)x_1^n;
  \end{array}
\right.$
\item[\rm2)]\, If $x_1\neq x_2$, then
$\left\{
  \begin{array}{l }
  a_n=\frac{a-x_2}{x_1-x_2}x_1^n-\frac{a-x_1}{x_1-x_2}x_2^n,\vspace{0.8ex}\\
  b_n=\frac{b}{x_1-x_2}(x_1^n-x_2^n),\vspace{0.8ex} \\
  c_n=\frac{c}{x_1-x_2}(x_1^n-x_2^n),\vspace{0.8ex} \\
  d_n=\frac{d-x_2}{x_1-x_2}x_1^n-\frac{d-x_1}{x_1-x_2}x_2^n.
  \end{array}
\right.$
\end{enumerate}
\end{lemma}

\begin{lemma}\label{l10}
Let $A=\begin{pmatrix} a & b \\ c & 0 \end{pmatrix}\in M_2\left(\mathbb{Z}\right)$ be a given matrix such that $bc\neq0$ and $\gcd(a,\, b,\,c)=1$. Let $K=\mathbb{Q}\left(\sqrt{a^2+4bc}\right)$ and let $\mathcal{O}_K$ be its ring of integers. If equation \eqref{e11} has a non-trivial solution in $C(A)$, then equation \eqref{e12} has a non-trivial solution in $\mathcal{O}_K$.
\end{lemma}
\begin{prof}
Suppose that $\left(X,\, Y,\, Z\right)$ is a non-trivial solution of equation \eqref{e11} in $C(A)$. By Proposition \ref{p6}, we obtain that $X,\, Y$ and $Z$ are pairwise commuting. Then there exists an invertible  matrix $P\in M_2(\mathbb{C})$ which simultaneously upper triangularizes the matrices $X,\, Y$ and $Z$. The assumption $uX^i+vY^j=wZ^k$ implies that $u\left(PXP^{-1}\right)^i+v\left(PYP^{-1}\right)^j=w\left(PZP^{-1}\right)^k$. We obtain
$$
u\begin{pmatrix} x_1 & \ast \\ 0 & x_2\end{pmatrix}^i+v\begin{pmatrix} y_1 & \ast \\ 0 & y_2\end{pmatrix}^j=w\begin{pmatrix} z_1 & \ast \\ 0 & z_2\end{pmatrix}^k.
$$
Then
$$
u\begin{pmatrix} x_1^i & \ast \\ 0 & x_2^i\end{pmatrix}+v\begin{pmatrix} y_1^j & \ast \\ 0 & y_2^j\end{pmatrix}=w\begin{pmatrix} z_1^k & \ast \\ 0 & z_2^k\end{pmatrix},
$$
where $x_s,\,y_s,\,z_s,\,s=1,\,2$ are the eigenvalues of $X,\, Y$ and $Z$, respectively. Comparing both sides, we have
$$
ux_s^i+vy_s^j=wz_s^k,\quad s=1,\,2.
$$
Therefore, $(x_s,\, y_s,\, z_s),\,s=1,\,2$ are non-trivial solutions of equation \eqref{e12} in $\mathcal{O}_K$.
\end{prof}

\begin{theorem}\label{th4}
Let $A=\begin{pmatrix} a & b \\ c & 0 \end{pmatrix}\in M_2\left(\mathbb{Z}\right)$ be a given matrix such that $bc\neq0$ and $\gcd(a,\, b,\, c)=1$. Let $K=\mathbb{Q}\left(\sqrt{a^2+4bc}\right)$ and let $\mathcal{O}_K$ be its ring of integers. Then the following statements hold.
\begin{enumerate}
\item[\rm1)]\, If $a^2+4bc$ is a square, then equation \eqref{e11} has a non-trivial solution in $C(A)$ if and only if equation \eqref{e12} has a non-trivial solution in $\mathbb{Z}$;
\item[\rm2)]\, If $a^2+4bc$ is not a square and $D$ is the unique square-free integer such that $a^2+4bc=m^2D$ for some $m\in\mathbb{N}$, then equation \eqref{e11} has a non-trivial solution in $C(A)$ if and only if equation \eqref{e12} has a non-trivial solution $\left(x,\, y,\, z\right)$ in $\mathcal{O}_K$ such that $x,\, y,\, z$ can be written in the form
$$
\frac{s+t\sqrt{D}}{2},\quad s,\, t\in\mathbb{Z},\, m\mid t.
$$
\end{enumerate}
\end{theorem}
\begin{prof}
1)\, In this case, we have $K=\mathbb{Q}\left(\sqrt{a^2+4bc}\right)=\mathbb{Q}$ and $\mathcal{O}_K=\mathbb{Z}$. Necessity follows from Lemma \ref{l10}. We next prove sufficiency. Assume that $\left(x,\, y,\, z\right)$ is a non-trivial solution of equation \eqref{e12} in $\mathbb{Z}$. Let $X=xI$, $Y=yI$ and $Z=zI$. Then $\left(X,\, Y,\, Z\right)$  is a non-trivial solution of equation \eqref{e11} in $C(A)$.

2)\, In this case, we have $K=\mathbb{Q}\left(\sqrt{a^2+4bc}\right)=\mathbb{Q}\left(\sqrt{D}\right)$. We now prove necessity. Assume that $\left(X,\, Y,\, Z\right)$ is a non-trivial solution of equation \eqref{e11} in $C(A)$. Then
\[
X=f_1I+g_1A,\, Y=f_2I+g_2A,\,Z=f_3I+g_3A
\]
for some $f_1,\, f_2,\, f_3,\,g_1,\, g_2,\, g_3\in \mathbb{Z}$. Let $x$, $y$ and $z$ be the eigenvalues of $X,\, Y$ and $Z$, respectively. Then
$$
x=\frac{tr\left(X\right)\pm g_1m\sqrt{D}}{2},\, y=\frac{tr\left(Y\right)\pm g_2m\sqrt{D}}{2},\, z=\frac{tr\left(Z\right)\pm g_3m\sqrt{D}}{2}.
$$
From Lemma \ref{l10}, it follows that $(x,\, y,\, z)$ is a non-trivial solution of equation \eqref{e12} in $\mathcal{O}_K$.

We next prove sufficiency. Assume that $\left(x_1,\, x_2,\, x_3\right)$ is a non-trivial solution of equation \eqref{e12} in $\mathcal{O}_K$ such that $x_1,\, x_2,\, x_3$ can be written in the form $\left(s+t\sqrt{D}\right)\big/2,\, s,\, t\in\mathbb{Z},\, m\mid t$. Let
$$x_r=\frac{s_r+t_r\sqrt{D}}{2},\quad s_r,\, t_r\in\mathbb{Z},\, m\mid t_r,\,r=1,\, 2,\, 3.$$
For $r=1,\, 2,\, 3$, let $\alpha_r=t_r/m$ and $\beta_r=\left(s_r-\alpha_ra\right)/2$. Since $m\mid t_r$, we have $\alpha_r\in\mathbb{Z}$ for $r=1,\,2,\, 3$. From $a^2+4bc=m^2D$, we obtain
\begin{equation}\label{e14}
\left(\alpha_ra\right)^2+4\alpha_r^2bc=t_r^2D,\quad r=1,\, 2,\, 3.
\end{equation}
If $D\equiv2,\, 3\pmod4$, then $2\mid s_r$ and $2\mid t_r$. By \eqref{e14}, we have $2\mid \alpha_ra$. Then $2\mid \left(s_r-\alpha_ra\right)$, i.e., $\beta_r\in\mathbb{Z}$. If $D\equiv1\pmod4$, then $2\mid \left(s_r+t_r\right)$. From \eqref{e14}, it follows that $\alpha_ra$ and $t_r$ have the same parity. Then $2\mid \left(s_r-\alpha_ra\right)$, i.e., $\beta_r\in\mathbb{Z}$. Hence, in any case, we have $\beta_r\in\mathbb{Z}$ for $r=1,\, 2,\, 3$. Let $$X_r=\beta_rI+\alpha_rA,\quad r=1,\, 2,\, 3.$$
By Lemma \ref{l5}, we have $X_r\in C(A)$ for $r=1,\, 2,\, 3$. We next show that $\left(X_1,\, X_2,\, X_3\right)$ is a non-trivial solution of equation \eqref{e11}. For $r=1,\, 2,\, 3$, notice that the eigenvalues of $X_r$ are $x_r$ and $\overline{x_r}$, where $\overline{x_r}$ denotes the conjugate of $x_r$. For a positive integer $n$, let
$$
X_r^n=\begin{pmatrix} a_{r, n} & b_{r, n} \\ c_{r, n} & d_{r, n}\end{pmatrix},\quad r=1,\, 2,\, 3.
$$
By Lemma \ref{le5}, we have
\begin{equation}\label{e16}
\left\{
\begin{array}{lr}
a_{r, n}=\frac{\left(a+m\sqrt{D}\right)x_r^n-\left(a-m\sqrt{D}\right)\overline{x_r}^n}{2m\sqrt{D}},\vspace{1ex} &  \\
b_{r, n}=\frac{b\left(x_r^n-\overline{x_r}^n\right)}{m\sqrt{D}},\vspace{1ex} &\\
c_{r, n}=\frac{c\left(x_r^n-\overline{x_r}^n\right)}{m\sqrt{D}},\vspace{1ex} &\\
d_{r, n}=\frac{\left(a+m\sqrt{D}\right)\overline{x_r}^n-\left(a-m\sqrt{D}\right)x_r^n}{2m\sqrt{D}} &\\
\end{array}
\right.
\end{equation}
for $r=1,\, 2,\, 3$ and $n\in\mathbb{N}$. Since $\left(x_1,\, x_2,\, x_3\right)$ is a non-trivial solution of equation \eqref{e12}, we have
\begin{equation}\label{e15}
ux_1^i+vx_2^j=wx_3^k.
\end{equation}
By \eqref{e16} and \eqref{e15}, we get
$$
\begin{aligned}
uX_1^i+vX_2^j&=u\begin{pmatrix} a_{1, i} & b_{1, i} \\ c_{1, i} & d_{1, i}\end{pmatrix}+v\begin{pmatrix} a_{2, j} & b_{2, j} \\ c_{2, j} & d_{2, j}\end{pmatrix}\\
&=\begin{pmatrix} ua_{1, i}+va_{2, j} & ub_{1, i}+vb_{2, j} \\ uc_{1, i}+vc_{2, j} & ud_{1, i}+vd_{2, j}\end{pmatrix}=\begin{pmatrix} wa_{3, k} & wb_{3, k}\\ wc_{3, k}& wd_{3, k}\end{pmatrix}=wX_3^k,
\end{aligned}
$$
which implies that $\left(X_1,\, X_2,\, X_3\right)$ is a non-trivial solution of equation \eqref{e11}.
\end{prof}

Let $i=j=k=n$. Then equations \eqref{e11} and \eqref{e12} become
\begin{equation}\label{e8}
uX^n+vY^n=wZ^n,\quad n\in\mathbb{N}
\end{equation}
and
\begin{equation}\label{e10}
ux^n+vy^n=wz^n,\quad n\in\mathbb{N},
\end{equation}
respectively.
\begin{theorem}\label{c5}
Let $A=\begin{pmatrix} a & b \\ c & 0 \end{pmatrix}\in M_2\left(\mathbb{Z}\right)$ be a given matrix such that $bc\neq0$ and $\gcd(a,\, b,\, c)=1$. Let $K=\mathbb{Q}\left(\sqrt{a^2+4bc}\right)$ and let $\mathcal{O}_K$ be its ring of integers. Then equation \eqref{e8} has a non-trivial solution in $C(A)$ if and only if equation \eqref{e10} has a non-trivial solution in $\mathcal{O}_K$.
\end{theorem}
\begin{prof}
If $a^2+4bc$ is a square, then the statement of theorem follows from Theorem \ref{th4} 1). Let us assume that $a^2+4bc$ is not a square. Let $D$ be the unique square-free integer such that
\begin{equation}\label{e3}
a^2+4bc=m^2D
\end{equation}
for some $m\in\mathbb{N}$. Necessity follows from Lemma \ref{l10}. We next prove sufficiency. Assume that $(x_1,\, x_2,\, x_3)$ is a non-trivial solution of equation \eqref{e10} in $\mathcal{O}_K$. Then $x_1,\, x_2,\, x_3$ can be written in the form $\left(s+t\sqrt{D}\right)\big/2,\, s,\, t\in\mathbb{Z}$. Let
$$x_i=\frac{s_i+t_i\sqrt{D}}{2},\quad s_i,\, t_i\in \mathbb{Z},\, i=1,\,2,\,3.$$
For $i=1,\, 2,\, 3$, let $\alpha_i=\left(ms_i-at_i\right)/2$. If $D\equiv2,\, 3\pmod4$, then $2\mid s_i$ and $2\mid t_i$. So $2\mid\left(ms_i-at_i\right)$, i.e., $\alpha_i\in\mathbb{Z}$. If $D\equiv1\pmod4$, then $2\mid \left(s_i+t_i\right)$. From \eqref{e3}, it follows that $a$ and $m$ have the same parity. Then $2\mid\left(ms_i-at_i\right)$, i.e., $\alpha_i\in\mathbb{Z}$. Hence, in any case, we have $\alpha_i\in\mathbb{Z}$ for $i=1,\, 2,\, 3$. Let
$$X_i=\alpha_iI+t_iA,\quad i=1,\, 2,\, 3.$$
By Lemma \ref{l5}, we have $X_i\in C(A)$ for $i=1,\, 2,\, 3$. We next show that $\left(X_1,\, X_2,\, X_3\right)$ is a non-trivial solution of equation \eqref{e8}. For $i=1,\, 2,\, 3$, notice that the eigenvalues of $X_i$ are $mx_i$ and $m\overline{x_i}$, where $\overline{x_i}$ denotes the conjugate of $x_i$.  For a positive integer $n$, let
$$
X_i^n=\begin{pmatrix} a_{i, n} & b_{i, n} \\ c_{i, n} & d_{i, n}\end{pmatrix},\quad i=1,\, 2,\, 3.
$$
By Lemma \ref{le5}, we have
\begin{equation}\label{e19}
\left\{
\begin{array}{lr}
a_{i, n}=\frac{\left(a+m\sqrt{D}\right)\left(mx_i\right)^n-\left(a-m\sqrt{D}\right)\left(m\overline{x_i}\right)^n}{2m\sqrt{D}},\vspace{1ex} &  \\
b_{i, n}=\frac{b\left(\left(mx_i\right)^n-\left(m\overline{x_i}\right)^n\right)}{m\sqrt{D}},\vspace{1ex} &\\
c_{i, n}=\frac{c\left(\left(mx_i\right)^n-\left(m\overline{x_i}\right)^n\right)}{m\sqrt{D}},\vspace{1ex} &\\
d_{i, n}=\frac{\left(a+m\sqrt{D}\right)\left(m\overline{x_i}\right)^n-\left(a-m\sqrt{D}\right)\left(mx_i\right)^n}{2m\sqrt{D}} &\\
\end{array}
\right.
\end{equation}
for $i=1,\, 2,\, 3$ and $n\in\mathbb{N}$. Since $\left(x_1,\, x_2,\, x_3\right)$ is a non-trivial solution of equation \eqref{e10}, we have
$$
ux_1^n+vx_2^n=wx_3^n,
$$
which implies that
\begin{equation}\label{e20}
u\left(mx_1\right)^n+v\left(mx_2\right)^n=w\left(mx_3\right)^n.
\end{equation}
By \eqref{e19} and \eqref{e20}, we get
$$
\begin{aligned}
uX_1^n+vX_2^n&=u\begin{pmatrix} a_{1, n} & b_{1, n} \\ c_{1, n} & d_{1, n}\end{pmatrix}+v\begin{pmatrix} a_{2, n} & b_{2, n} \\ c_{2, n} & d_{2, n}\end{pmatrix}\\
&=\begin{pmatrix} ua_{1, n}+va_{2, n} & ub_{1, n}+vb_{2, n} \\ uc_{1, n}+vc_{2, n} & ud_{1, n}+vd_{2, n}\end{pmatrix}=\begin{pmatrix} wa_{3, n} & wb_{3, n}\\ wc_{3, n}& wd_{3, n}\end{pmatrix}=wX_3^n,
\end{aligned}
$$
which implies that $\left(X_1,\, X_2,\, X_3\right)$ is a non-trivial solution of equation \eqref{e8}.
\end{prof}

Let $u=v=w=1$ and $n\geq3$. Then equations \eqref{e8} and \eqref{e10} become the Fermat's matrix equation
\begin{equation}\label{e17}
X^n+Y^n=Z^n,\quad n\in\mathbb{N},\, n\geq3
\end{equation}
and the Fermat's equation
\begin{equation}\label{e18}
x^n+y^n=z^n,\quad n\in\mathbb{N},\, n\geq3,
\end{equation}
respectively.
\begin{corollary}\label{c1}
Let $A=\begin{pmatrix} a & b \\ c & 0 \end{pmatrix}\in M_2\left(\mathbb{Z}\right)$ be a given matrix such that $bc\neq0$ and $\gcd(a,\, b,\, c)=1$. Let $K=\mathbb{Q}\left(\sqrt{a^2+4bc}\right)$ and let $\mathcal{O}_K$ be its ring of integers. Then equation \eqref{e17} has a non-trivial solution in $C(A)$ if and only if equation \eqref{e18} has a non-trivial solution in $\mathcal{O}_K$.
\end{corollary}
\begin{prof}
Directly from Theorem \ref{c5}.
\end{prof}

From Corollary \ref{c1}, we conclude that the solvability of the Fermat's matrix equation \eqref{e17} in $C(A)$ can be reduced to the solvability of the Fermat's equation \eqref{e18} in quadratic fields.  However, the  solvability of the Fermat's equation in quadratic fields is unsolved. The following lemmas list some known results about the solvability of the Fermat's equation  in quadratic fields.

\begin{lemma}{\rm (\cite{F.})}\label{l1}
Equation \eqref{e18} has no non-trivial solutions in $\mathbb{Q}\left(\sqrt{2}\right)$ for $n\geq4$.
\end{lemma}

\begin{lemma}{\rm (\cite{Freitas})}\label{l2}
Let $3\leq D\neq5,\, 17\leq23$ be a square-free integer. Then equation \eqref{e18} has no non-trivial solutions in $\mathbb{Q}\left(\sqrt{D}\right)$ for $n\geq4$.
\end{lemma}

\begin{lemma}{\rm (\cite{A.Aigner})}\label{l3}
Let $D\neq1$ be a square-free integer. Then the equation $x^4+y^4=z^4$ has non-trivial solutions in $\mathbb{Q}\left(\sqrt{D}\right)$ if and only if $D=-7$, and all non-trivial solutions in $\mathbb{Q}\left(\sqrt{-7}\right)$ can be reduced to the  solution
$$\left(\frac{1+\sqrt{-7}}{2}\right)^4+\left(\frac{1-\sqrt{-7}}{2}\right)^4=1.$$
\end{lemma}
\begin{lemma}{\rm (\cite{Die})}\label{l4}
Equation \eqref{e18} has no non-trivial solutions in quadratic fields for $n=6,\, 9$.
\end{lemma}

Combining the above results, we have the following corollaries.
\begin{corollary}\label{c3}
Let $2\leq D\neq5,\, 17\leq23$ be a square-free integer. Let $a,\, b,\, c,\, m$ be integers such that $a^2+4bc=m^2D$, $bc\neq0$ and $\gcd\left(a,\, b,\, c\right)=1$. Then equation \eqref{e17} has no non-trivial solutions in $C(A)$ for $n\geq4$, where $A=\begin{pmatrix} a & b \\ c & 0\end{pmatrix}$.
\end{corollary}
\begin{prof}
Let $K=\mathbb{Q}\left(\sqrt{a^2+4bc}\right)$ and let $\mathcal{O}_K$ be its ring of integers. By Corollary \ref{c1}, equation \eqref{e17} has  a non-trivial solution in $C(A)$ if and only if equation \eqref{e18} has a non-trivial solution in $\mathcal{O}_K$. If $m=0$, then $K=\mathbb{Q}$ and $\mathcal{O}_K=\mathbb{Z}$. By Fermat's last theorem, equation \eqref{e18} has no non-trivial solutions in $\mathcal{O}_K$. Then equation \eqref{e17} has no non-trivial solutions in $C(A)$. If $m\neq0$, then $K=\mathbb{Q}\left(\sqrt{D}\right)$. By Lemmas \ref{l1} and \ref{l2}, equation \eqref{e18} has no non-trivial solutions in $\mathcal{O}_K$ for $2\leq D\neq5,\, 17\leq23$ when $n\geq4$. Then equation \eqref{e17} has no non-trivial solutions in $C(A)$ for $n\geq4$.
\end{prof}

\begin{example}\label{ex1}
\rm
Let $2\leq D\neq5,\, 17\leq23$ be a square-free integer. Let $a,\,m$ be integers such that $a^2+4=m^2D$. From Corollary \ref{c3}, it follows that equation \eqref{e17} has no  non-trivial solutions in $C(A)$ for $n\geq4$, where $A=\begin{pmatrix} a & 1 \\ 1 & 0\end{pmatrix}$. For example, let
$$\left(D,\,a,\,m\right)=\left(2,\,\pm2,\,2\right),\,\left(13,\,\pm3,\,1\right),\,\left(10,\,\pm6,\,2\right).$$
Then equation \eqref{e17} has no  non-trivial solutions in $C(A)$ for $n\geq4$, where $A=\begin{pmatrix} a & 1 \\ 1 & 0\end{pmatrix}$ and $a=\pm2,\, \pm3,\,\pm6$.
\end{example}

\begin{corollary}\label{c2}
Let $A=\begin{pmatrix} a & b \\ c & 0 \end{pmatrix}\in M_2(\mathbb{Z})$ be a given matrix such that $bc\neq0$ and $\gcd(a,\, b,\, c)=1$. Then the following statements hold.
\begin{enumerate}
\item[\rm1)]\, If $a^2+4bc$ is a square, then equation \eqref{e17} has no non-trivial solutions in $C(A)$;
\item[\rm2)]\, The equation $X^4+Y^4=Z^4$ has a non-trivial solution in $C(A)$ if and only if $$\mathbb{Q}\left(\sqrt{a^2+4bc}\right)=\mathbb{Q}\left(\sqrt{-7}\right);$$
\item[\rm3)]\, Equation \eqref{e17} has no non-trivial solutions in $C(A)$ for $n=6,\, 9$;
\item[\rm4)]\, If equation \eqref{e17} has at least one non-trivial solution in $C(A)$, then it has infinitely many non-trivial solutions in $C(A)$.
\end{enumerate}
\end{corollary}
\begin{prof}
Let $K=\mathbb{Q}\left(\sqrt{a^2+4bc}\right)$ and let $\mathcal{O}_K$ be its ring of integers. By Corollary \ref{c1}, equation \eqref{e17} has  a non-trivial solution in $C(A)$ if and only if equation \eqref{e18} has a non-trivial solution in $\mathcal{O}_K$.

1)\, In this case, we have $K=\mathbb{Q}$ and $\mathcal{O}_K=\mathbb{Z}$. By Fermat's last theorem, equation \eqref{e18} has no non-trivial solutions in $\mathcal{O}_K$. Then equation \eqref{e17} has no non-trivial solutions in $C(A)$.

2)\, By Lemma \ref{l3}, the equation $x^4+y^4=z^4$ has a non-trivial solution in $\mathcal{O}_K$ if and only if $K=\mathbb{Q}\left(\sqrt{-7}\right)$. Therefore, the equation $X^4+Y^4=Z^4$ has a non-trivial solution in $C(A)$ if and only if $K=\mathbb{Q}\left(\sqrt{-7}\right)$.

3)\, By Lemma \ref{l4}, equation \eqref{e18} has no non-trivial solutions in $\mathcal{O}_K$ for $n=6,\, 9$. Then equation \eqref{e17} has no non-trivial solutions in $C(A)$ for $n=6,\, 9$.

4)\, Suppose that $\left(X,\, Y,\, Z\right)$ is a non-trivial solution of equation \eqref{e17} in $C(A)$. From 1), it follows that $a^2+4bc$ is not a square. By Corollary \ref{c0}, we know that $C(A)$ forms an integral domain under the operations of matrix addition and multiplication. Let $B\in C(A)$ be an arbitrary matrix such that $B\neq O$. Then by the proof of Proposition \ref{p1}, we obtain $\det(B)\neq0$. Since $\left(X,\, Y,\, Z\right)$ is a non-trivial solution of equation \eqref{e17} in $C(A)$, we have
$$
\left(BX\right)^n+\left(BY\right)^n=B^nX^n+B^nY^n=B^n\left(X^n+Y^n\right)=B^nZ^n=\left(BZ\right)^n.
$$
This means that $\left(BX,\, BY,\, BZ\right)$ are non-trivial solutions of equation \eqref{e17} in $C(A)$. Since $C(A)$ has no zero divisors, these non-trivial solutions are pairwise different.
\end{prof}

\begin{example}\label{ex2}
\rm Let $q$ be an integer and let $A=\begin{pmatrix} q & 1 \\ 1 & 0 \end{pmatrix}$. Notice that $\mathbb{Q}\left(\sqrt{q^2+4}\right)\neq\mathbb{Q}\left(\sqrt{-7}\right)$. From Corollary \ref{c2} 2), it follows that equation \eqref{e17} has no  non-trivial solutions in $C(A)$ for $n=4$. Moreover, by Corollary \ref{c2} 3), we know that equation \eqref{e17} has no  non-trivial solutions in $C(A)$ for $n=6,\, 9$. Therefore, equation \eqref{e17} has no  non-trivial solutions in $C(A)$ for $n=4,\, 6,\, 9$.
\end{example}
\begin{remark}
\rm Examples \ref{ex1} and \ref{ex2} are given in \cite[Theorems 3 and 5]{ch}.
\end{remark}

Let $K$ be a quadratic field and $\mathcal{O}_K$ its ring of integers. Let $D$ be the unique square-free integer such that $K=\mathbb{Q}\left(\sqrt{D}\right)$. From a given non-trivial solution of equation \eqref{e18} in $\mathcal{O}_K$, we can construct infinitely many classes of $2\times 2$ matrices such that equation \eqref{e17} has non-trivial solutions in these classes. Assume that
$$
\left(\frac{s_1+t_1\sqrt{D}}{2}\right)^n+\left(\frac{s_2+t_2\sqrt{D}}{2}\right)^n=\left(\frac{s_3+t_3\sqrt{D}}{2}\right)^n
$$
is a given non-trivial solution of equation \eqref{e18} in $\mathcal{O}_K$, where $s_i,\, t_i\in\mathbb{Z},\, i=1,\, 2,\, 3.$ Let $a,\, b,\, c$ be integers and $m$ a positive integer such that $a^2+4bc=m^2D$, $bc\neq0$ and $\gcd\left(a,\, b,\, c\right)=1$. Indeed, there are infinitely many such $a,\, b,\, c,\, m$. Let $t$ be an arbitrary positive integer. If $D\equiv1\pmod4$ and $D=1+4k$ for some $k\in\mathbb{Z}$, then $\left(a,\, b,\, c,\, m\right)=\left(t,\, 1,\, kt^2,\, t\right)$ satisfy the above conditions. If $D\equiv2\pmod4$ and $D=2+4k$ for some $k\in\mathbb{Z}$, then $\left(a,\, b,\, c,\, m\right)=\left(2t,\, 1,\, t^2(1+4k),\, 2t\right)$ satisfy the above conditions. If $D\equiv3\pmod4$ and $D=3+4k$ for some $k\in\mathbb{Z}$, then $\left(a,\, b,\, c,\, m\right)=\left(2t,\, 1,\, 2t^2(1+2k),\, 2t\right)$ satisfy the above conditions. From the proof of Theorem \ref{c5}, it follows that
$$
\begin{pmatrix} \frac{ms_1+at_1}{2} & t_1b \\ t_1c & \frac{ms_1-at_1}{2} \end{pmatrix}^n+\begin{pmatrix} \frac{ms_2+at_2}{2} & t_2b \\ t_2c & \frac{ms_2-at_2}{2} \end{pmatrix}^n=\begin{pmatrix} \frac{ms_3+at_3}{2} & t_3b \\ t_3c & \frac{ms_3-at_3}{2} \end{pmatrix}^n
$$
are non-trivial solutions of equation \eqref{e17}, and the corresponding matrix classes are $C(A)$, where $A=\begin{pmatrix} a & b \\ c & 0\end{pmatrix}$. Next, we give some examples to illustrate how to construct non-trivial solutions of the Fermat's matrix equation in $M_2\left(\mathbb{Z}\right)$ from a given equality in this manner. Moreover, we have not found other similar methods for constructing non-trivial solutions in the literature.
\begin{example}\label{eq7}
\rm
In \cite[Theorem 2]{ch}, M. T. Chien and J. Meng gave a non-trivial solution of the equation $X^3+Y^3=Z^3$ in $M_2\left(\mathbb{Z}\right)$:
\begin{equation}\label{eq4}
\begin{pmatrix} 7 & 3 \\ 3 & 4 \end{pmatrix}^3+\begin{pmatrix} 11 & 6 \\ 6 & 5 \end{pmatrix}^3=\begin{pmatrix} 12 & 6 \\ 6 & 6 \end{pmatrix}^3.
\end{equation}
Note that their eigenvalues satisfy the equality
$$
\left(\frac{11+3\sqrt{5}}{2}\right)^3+\left(8+3\sqrt{5}\right)^3=\left(9+3\sqrt{5}\right)^3.
$$
From this equality, we can construct infinitely many classes of $2\times 2$ matrices such that the equation $X^3+Y^3=Z^3$ has non-trivial solutions in these classes. Let $a,\, b,\, c$ be integers and $m$ a positive integer such that $a^2+4bc=5m^2$, $bc\neq0$ and $\gcd\left(a,\, b,\, c\right)=1$. Then
\begin{equation}\label{eq6}
\begin{pmatrix} \frac{11m+3a}{2} & 3b \\ 3c & \frac{11m-3a}{2} \end{pmatrix}^3+\begin{pmatrix} 8m+3a & 6b \\ 6c & 8m-3a \end{pmatrix}^3=\begin{pmatrix} 9m+3a & 6b \\ 6c & 9m-3a \end{pmatrix}^3
\end{equation}
are non-trivial solutions of the equation $X^3+Y^3=Z^3$, and the corresponding matrix classes are $C(A)$, where $A=\begin{pmatrix} a & b \\ c & 0\end{pmatrix}$. Let $a=b=c=m=1$. Then we get the non-trivial solution \eqref{eq4}, and the corresponding matrix class is $C(B)$, where $B=\begin{pmatrix} 1 & 1 \\ 1 & 0\end{pmatrix}$.
\end{example}

Example \ref{eq7} shows that the non-trivial solution \eqref{eq4} can be obtained from \eqref{eq6} and there are infinitely many such non-trivial solutions.
\begin{example}\label{ex3}
\rm
In \cite{Bu}, W. Burnside gave the equality
\begin{equation}\label{e22}
\left(-3+\sqrt{-3\left(1+4k^3\right)}\right)^3+\left(-3-\sqrt{-3\left(1+4k^3\right)}\right)^3=\left(6k\right)^3,
\end{equation}
where $k\neq0,\, -1$ is an integer. We claim that $-3\left(1+4k^3\right)$ is not a square. Otherwise, $-3\left(1+4k^3\right)=q^2$ for some $q\in\mathbb{N}$. From \eqref{e22}, it follows that
$$
\left(-3+q\right)^3+\left(-3-q\right)^3=\left(6k\right)^3.
$$
By Fermat's last theorem, we have $q=3$. This implies that $-3\left(1+4k^3\right)=9$, so we obtain $k=-1$, a contradiction to $k\neq0,\, -1$. Hence, $-3\left(1+4k^3\right)$ is not a square. Let $D\neq1$ be the unique square-free integer such that $-3\left(1+4k^3\right)=t^2D$ for some $t\in\mathbb{N}$. Then equation \eqref{e22} becomes
$$
\left(-3+t\sqrt{D}\right)^3+\left(-3-t\sqrt{D}\right)^3=\left(6k\right)^3.
$$
From this equality, we can construct infinitely many classes of $2\times 2$ matrices such that the equation $X^3+Y^3=Z^3$ has non-trivial solutions in these classes. Let $a,\, b,\, c$ be integers and $m$ a positive integer such that $a^2+4bc=m^2D$, $bc\neq0$ and $\gcd\left(a,\, b,\, c\right)=1$. Then
$$
\begin{pmatrix} -3m+at & 2tb \\ 2tc & -3m-at \end{pmatrix}^3+\begin{pmatrix} -3m-at & -2tb \\ -2tc & -3m+at \end{pmatrix}^3=\begin{pmatrix} 6mk & 0 \\ 0 & 6mk \end{pmatrix}^3
$$
are non-trivial solutions of the equation $X^3+Y^3=Z^3$, and the corresponding matrix classes are $C(A)$, where $A=\begin{pmatrix} a & b \\ c & 0\end{pmatrix}$.
\end{example}

\begin{example}\label{ex4}
\rm
In \cite{A.Aigner}, A. Aigner gave the equality
\begin{equation}\label{eq8}
\left(\frac{1+\sqrt{-7}}{2}\right)^4+\left(\frac{1-\sqrt{-7}}{2}\right)^4=1.
\end{equation}
From this equality, we can construct infinitely many classes of $2\times 2$ matrices such that the equation $X^4+Y^4=Z^4$ has non-trivial solutions in these classes.  Let $a,\, b,\, c$ be integers and $m$ a positive integer such that $a^2+4bc=-7m^2$, $bc\neq0$ and $\gcd\left(a,\, b,\, c\right)=1$. Then
$$
\begin{pmatrix} \frac{m+a}{2} & b \\ c & \frac{m-a}{2} \end{pmatrix}^4+\begin{pmatrix} \frac{m-a}{2} & -b \\ -c & \frac{m+a}{2} \end{pmatrix}^4=\begin{pmatrix} m & 0 \\ 0 & m \end{pmatrix}^4
$$
are non-trivial solutions of the equation $X^4+Y^4=Z^4$, and the corresponding matrix classes are $C(A)$, where $A=\begin{pmatrix} a & b \\ c & 0\end{pmatrix}$.
\end{example}

Examples \ref{ex3} and \ref{ex4} show that we can construct infinitely many non-trivial solutions of the Fermat's matrix equation with exponents $3$ and $4$ in $M_2\left(\mathbb{Z}\right)$ from the equalities \eqref{e22} and \eqref{eq8}, respectively.

\begin{example}\label{ex5}
\rm
Let $r$ and $s$ be arbitrary integers such that they are not all equal to zero. Let $k$ be an arbitrary positive integer. In \cite[Theorem 3]{km}, I. Kaddoura and B. Mourad proved that
\begin{equation}\label{eq3}
\begin{pmatrix} s & -r \\ r & s-r \end{pmatrix}^{n}+\begin{pmatrix} r-s & s \\ -s & r \end{pmatrix}^n=\begin{pmatrix} r & s-r \\ r-s & s \end{pmatrix}^n
\end{equation}
are non-trivial solutions of the equations $X^n+Y^n=Z^n$ in $M_2\left(\mathbb{Z}\right)$, where $n=6k+1,\,6k+5$. Next, we show that the non-trivial solutions \eqref{eq3} can be obtained from two equalities. For polynomials with integer coefficients, we have the following congruences \cite[Lemma 2]{km}.
\begin{equation}\label{eq1}
\left(x+y\right)^{6k+1}-x^{6k+1}-y^{6k+1}\equiv0\mod\left(xy+x^2+y^2\right)
\end{equation}
\begin{equation}\label{eq2}
\left(x+y\right)^{6k+5}-x^{6k+5}-y^{6k+5}\equiv0\mod\left(xy+x^2+y^2\right)^2
\end{equation}
Let $f(x,y)=xy+x^2+y^2$. Then $f\left((2s-r+r\sqrt{-3})/2,(2r-s-s\sqrt{-3})/2\right)=0$. From \eqref{eq1} and \eqref{eq2}, it follows that
$$
\left(\frac{2s-r+r\sqrt{-3}}{2}\right)^n+\left(\frac{2r-s-s\sqrt{-3}}{2}\right)^n=\left(\frac{r+s+(r-s)\sqrt{-3}}{2}\right)^n,
$$
where $n=6k+1,\,6k+5$. From these two equalities, we can construct infinitely many classes of $2\times 2$ matrices such that the equations $X^n+Y^n=Z^n,\,n=6k+1,\,6k+5$ have infinitely many non-trivial solutions in these classes. Let $a,\, b,\, c$ be integers and $m$ a positive integer such that $a^2+4bc=-3m^2$, $bc\neq0$ and $\gcd\left(a,\, b,\, c\right)=1$. Then
\begin{equation}\label{eq9}
\begin{pmatrix} \frac{m(2s-r)+ar}{2} & rb \\ rc & \frac{m(2s-r)-ar}{2} \end{pmatrix}^n+\begin{pmatrix} \frac{m(2r-s)-as}{2} & -sb \\ -sc & \frac{m(2r-s)+as}{2} \end{pmatrix}^n=\begin{pmatrix} \frac{m(r+s)+a(r-s)}{2} & (r-s)b \\ (r-s)c & \frac{m(r+s)-a(r-s)}{2} \end{pmatrix}^n
\end{equation}
are non-trivial solutions of the equations $X^n+Y^n=Z^n$ in $C(A)$, where $A=\begin{pmatrix} a & b \\ c & 0\end{pmatrix}$ and $n=6k+1,\,6k+5$. Let $a=c=m=1$ and $b=-1$. Then we get the non-trivial solutions \eqref{eq3}, and the corresponding matrix class is $C(B)$, where $B=\begin{pmatrix} 1 & -1 \\ 1 & 0\end{pmatrix}$.
\end{example}

Example \ref{ex5} shows that the non-trivial solutions \eqref{eq3} can be obtained from \eqref{eq9} and there are infinitely many such non-trivial solutions.

\section{Catalan's equation over $M_2(\mathbb{Z})$}\label{s2}
In this section, we study the solvability of the Catalan's matrix equation
\begin{equation}\label{e7}
X^m-Y^n=I,\quad m,\, n\in\mathbb{N},\, m,\,n\geq3
\end{equation}
in  $M_2\left(\mathbb{Z}\right)$. Here we require $m,\,n\geq3$. Indeed, if $m=2$ or $n=2$, without loss of generality, we can assume that $m=2$. For any integer $t\neq0,\, -1$, let $A$ be a matrix in $M_2\left(\mathbb{Z}\right)$ such that $tr(A)=0$ and $\det(A)=-t^n-1$. Then $A^2=\left(t^n+1\right)I$, i.e.,
$$
A^2-\left(tI\right)^n=I.
$$
Therefore, we can get the non-trivial solutions $\left(X,\, Y,\, m,\, n\right)=\left(A,\, tI,\, 2,\,n\right)$ of the Catalan's matrix equation. However, these solutions are trivial. Thus we assume that $m,\, n\geq3$.

\begin{definition}
Let $K$ be a quadratic field and $\mathcal{O}_K$ its ring of integers. Let $x\in \mathcal{O}_K$. If there is a positive integer $t$ such that $x^t\in\mathbb{Z}$, then we say that $x$ has finite exponent $t_0$, where $t_0$ is the smallest positive integer with such property. Otherwise, we say that $x$ has infinite exponent $\infty$. We denote the exponent of $x$  by $E(x)$.
\end{definition}

About exponent, we have the following statements.
\begin{proposition}
Let $K$ be a quadratic field and $\mathcal{O}_K$ its ring of integers. For $x\in \mathcal{O}_K$, we have
$$E(x)\in\{1,\, 2,\, 3,\, 4,\, 6,\, \infty\}.$$
\end{proposition}
\begin{prof}
Suppose that $x$ has finite exponent. If $x$ is an integer, then $x$ has exponent $1$. If $x$ is not an integer, then $x^{E(x)}\in\mathbb{Z}$ and $x^i\notin\mathbb{Z}$ for $1\leq i<E(x)$, which imply that $\overline{x}/x$ is a primitive $E(x)$th root of unity, where $\overline{x}$ denotes  the conjugate of $x$. We know that the degree of $\overline{x}/x$ over $\mathbb{Q}$ is $\varphi\left(E(x)\right)$, where $\varphi$ is Euler's totient function. Then $\varphi\left(E(x)\right)\leq2$, which implies that $E(x)\in\{2,\, 3,\, 4,\, 6\}$.
\end{prof}
\begin{proposition}\label{p3}
Let $K$ be a quadratic field and $\mathcal{O}_K$ its ring of integers. If $x\in \mathcal{O}_K$ has finite exponent, then for $n\in\mathbb{N}$, $x^n\in\mathbb{Z}$ if and only if $E(x)\mid n$.
\end{proposition}
\begin{prof}
The sufficiency is clear. We next prove necessity. The case $x\in\mathbb{Z}$ is evident, so we assume that $x\notin\mathbb{Z}$. Let $n=E(x)q+r$, where $q,\, r\in\mathbb{Z}$ and $0\leq r<E(x)$. Then we have $x^r=x^{n-E(x)q}$ and $\overline{x}^r=\overline{x}^{n-E(x)q}$, where $\overline{x}$ denotes  the conjugate of $x$. Since $x^n,\, x^{E(x)}\in\mathbb{Z}$, we botain
$$
\frac{\overline{x}^r}{x^r}=\frac{\overline{x}^{n-E(x)q}}{x^{n-E(x)q}}=\frac{\overline{x}^n}{x^n}\cdot\left(\frac{\overline{x}^{E(x)}}{x^{E(x)}}\right)^{-q}=1,
$$
which implies that $x^r\in\mathbb{Z}$. If $r\neq0$, then we obtain a contradiction to the minimality of $E(x)$. So $r=0$, which means that $E(x)\mid n$.
\end{prof}
\begin{proposition}\label{p4}
Let $K$ be a quadratic field and $\mathcal{O}_K$ its ring of integers. Let $D$ be the unique square-free integer such that $K=\mathbb{Q}\left(\sqrt{D}\right)$. Let $E_j=\left\{x\in \mathcal{O}_K: E(x)=j\right\}$, $j=1,\, 2,\, 3,\, 4,\, 6$ and let $i=\sqrt{-1}$, $\omega=\left(-1+\sqrt{-3}\right)/2$. Then the following statements hold.
\begin{enumerate}
\item[\rm1)]\, $E_1=\mathbb{Z}$;
\item[\rm2)]\, $E_2=\left\{k\sqrt{D}: k\in\mathbb{Z},\, k\neq0\right\}$;
\item[\rm3)]\, $E_3\neq\emptyset$ if and only if $D=-3$, and $E_3=\left\{k\omega,\, k\overline{\omega}: k\in\mathbb{Z},\, k\neq0\right\}$;
\item[\rm4)]\, $E_4\neq\emptyset$ if and only if $D=-1$, and $E_4=\left\{k(1+i),\, k(1-i): k\in\mathbb{Z},\, k\neq0\right\}$;
\item[\rm5)]\, $E_6\neq\emptyset$ if and only if $D=-3$, and $E_6=\left\{k(1-\omega),\, k(1-\overline{\omega}): k\in\mathbb{Z},\, k\neq0\right\}$.
\end{enumerate}
\end{proposition}
\begin{prof}
1)\, Clearly.

2)\, If $E(x)=2$, then $\overline{x}/x$ is a primitive $2$th root of unity, i.e., $\overline{x}/x=-1$. This implies that $x=k\sqrt{D}$, where $k$ is a nonzero integer.

3)\, If $E(x)=3$, then $\overline{x}/x$ is a primitive $3$th root of unity, i.e., $\overline{x}/x=\omega$ or $\overline{\omega}$. This implies that $x=k\omega$ or $k\overline{\omega}$, where $k$ is a nonzero integer.

4)\, If $E(x)=4$, then $\overline{x}/x$ is a primitive $4$th root of unity, i.e., $\overline{x}/x=\pm i$. This implies that $x=k(1+i)$ or $k(1-i)$, where $k$ is a nonzero integer.

5)\, If $E(x)=6$, then $\overline{x}/x$ is a primitive $6$th root of unity, i.e., $\overline{x}/x=-\omega$ or $-\overline{\omega}$. This implies that $x=k(1-\omega)$ or $k\left(1-\overline{\omega}\right)$, where $k$ is a nonzero integer.
\end{prof}

\begin{lemma}\label{l7}
If equation \eqref{e7} has a non-trivial solution in $M_2\left(\mathbb{Z}\right)$, then the equation
\begin{equation}\label{e9}
x^m-y^n=1,\quad m,\, n\in\mathbb{N},\, m,\,n\geq3
\end{equation}
has a non-trivial solution in algebraic integers $x$ and $y$ of degree less than or equal to $2$.
\end{lemma}
\begin{prof}
Suppose that $\left(X,\, Y,\, m,\, n\right)$ is a non-trivial solution of equation \eqref{e7} in $M_2\left(\mathbb{Z}\right)$. There exists an invertible  matrix $P\in M_2(\mathbb{C})$ which  upper triangularizes the matrix $X$. The assumption $X^m-Y^n=I$ implies that $\left(PXP^{-1}\right)^m-\left(PYP^{-1}\right)^n=I$. We obtain
$$
\begin{pmatrix} x_1 & \ast \\ 0 & x_2\end{pmatrix}^m-\left(PYP^{-1}\right)^n=I.
$$
Then
\begin{equation}\label{e0}
\left(PYP^{-1}\right)^n=\begin{pmatrix} x_1^m-1 & \ast \\ 0 & x_2^m-1\end{pmatrix},
\end{equation}
where $x_s,\,s=1,\,2$ are the eigenvalues of $X$. Let $y_s,\,s=1,\,2$ be the eigenvalues of $Y$. Then the eigenvalues of $\left(PYP^{-1}\right)^n$ are $y_s^n,\,s=1,\,2$. By \eqref{e0}, we have
$$
y_s^n=x_s^m-1,\quad s=1,\,2.
$$
Therefore, $(x_s,\, y_s,\, m,\, n),\,s=1,\,2$ are non-trivial solutions of equation \eqref{e9}.
\end{prof}

Lemma \ref{l7} tells us that we should consider the solvability of equation \eqref{e9} in algebraic integers $x$ and $y$ of degree less than or equal to $2$.

\begin{lemma}\label{l8}
If $x$ or $y$ is an integer, then all non-trivial solutions of equation \eqref{e9} are
$$\left(\pm\sqrt{\pm3},\, 2,\, 4,\, 3\right).$$
\end{lemma}
\begin{prof}
We consider the following two cases.

{\bf Case 1:} If $x$ is an integer, then $y^n\in\mathbb{Z}$. So $y$  has finite exponent. By Proposition \ref{p3}, we obtain $E(y)\mid n$. By Proposition \ref{p4}, equation \eqref{e9} becomes
$$
\begin{cases}
 x^m-y^n=1, & \mbox{if}\, \, E(y)=1, \\
 x^m-\left(k^2D\right)^{n/2}=1, &  \mbox{if}\, \, E(y)=2, \\
 x^m-k^n=1, & \mbox{if} \,\, E(y)=3,\\
 x^m+(-1)^{n/4+1}\left(2k^2\right)^{n/2}=1, & \mbox{if} \,\, E(y)=4,\\
 x^m+(-1)^{n/6+1}\left(3k^2\right)^{n/2}=1, & \mbox{if} \,\, E(y)=6,\\
\end{cases}
$$
where $D\neq1$ is a square-free integer and $k$ is a nonzero integer. By Catalan's conjecture, we know that these equations have no non-trivial solutions.

{\bf Case 2:} If $y$ is an integer, then $x^m\in\mathbb{Z}$.  So $x$  has finite exponent. By Proposition \ref{p3}, we obtain $E(x)\mid m$. By Proposition \ref{p4}, equation \eqref{e9} becomes
$$
\begin{cases}
 x^m-y^n=1, & \mbox{if}\, \, E(x)=1, \\
 \left(k^2D\right)^{m/2}-y^n=1, &  \mbox{if}\, \, E(x)=2, \\
 k^m-y^n=1, & \mbox{if} \,\, E(x)=3,\\
 (-1)^{m/4}\left(2k^2\right)^{m/2}-y^n=1, & \mbox{if} \,\, E(x)=4,\\
 (-1)^{m/6}\left(3k^2\right)^{m/2}-y^n=1, & \mbox{if} \,\, E(x)=6,\\
\end{cases}
$$
where $D\neq1$ is a square-free integer and $k$ is a nonzero integer. By Catalan's conjecture, we know that only one of these equations has non-trivial solutions, and this equation is
$$\left(k^2D\right)^{m/2}-y^n=1.$$
Then all non-trivial solutions of this equation are $k^2D=\pm3,\, y=2,\, m/2=2,\, n=3$. So we obtain $\left(x,\, y,\, m,\, n\right)=\left(\pm\sqrt{\pm3},\, 2,\, 4,\, 3\right).$
\end{prof}
\begin{theorem}\label{th5}
If the eigenvalues of $X$ or $Y$ are integers, then all non-trivial solutions of equation \eqref{e7} in $M_2(\mathbb{Z})$ are given by $$\left\{\left(X,\, 2I,\, 4,\, 3\right): X\in M_2(\mathbb{Z}),\, tr(X)=0,\, \det(X)=\pm3\right\}.$$
\end{theorem}
\begin{prof}
Assume that $x_s,\,y_s,\,s=1,\,2$ are the eigenvalues of $X$ and $Y$, respectively. By Lemma \ref{l7}, we know that $(x_s,\, y_s,\, m,\, n),\,s=1,\,2$ are non-trivial solutions of equation \eqref{e9}. In this case, by Lemma \ref{l8}, we have
$$\left(x_s,\, y_s,\, m,\, n\right)\in\left\{\left(\pm\sqrt{\pm3},\, 2,\, 4,\, 3\right)\right\},\quad s=1,\,2.$$
By Lemma \ref{le5} and a direct computation, we can get all non-trivial solutions of equation \eqref{e7} in $M_2(\mathbb{Z})$ in this case, which are given in the theorem.
\end{prof}

\begin{lemma}\label{l9}
If $x$ and $y$ are quadratic algebraic integers such that $x^m$ is an integer, then all non-trivial solutions of equation \eqref{e9} are $$\left(\pm\sqrt{\pm3},\, \pm\sqrt{2},\, 4,\, 6\right),\, \left(\pm\sqrt{\pm3},\, 2\omega,\, 4,\, 3\right),\, \left(\pm\sqrt{\pm3},\, 2\overline{\omega},\, 4,\, 3\right),$$ where $\omega=(-1+\sqrt{-3})/2$.
\end{lemma}
\begin{prof}
Since $x^m$ is an integer, it follows that $y^n$ is also an integer. So $x$ and $y$  have finite exponent. By Proposition \ref{p3}, we have $E(x)\mid m$ and $E(y)\mid n$. We next consider the following four cases.

{\bf Case 1:} $E(x)=2$.  By Proposition \ref{p4}, equation \eqref{e9} becomes
$$
\begin{cases}
 \left(k_1^2D_1\right)^{m/2}-\left(k_2^2D_2\right)^{n/2}=1, &  \mbox{if}\, \, E(y)=2, \\
 \left(k_1^2D_1\right)^{m/2}-k_2^n=1, & \mbox{if} \,\, E(y)=3,\\
 \left(k_1^2D_1\right)^{m/2}+(-1)^{n/4+1}\left(2k_2^2\right)^{n/2}=1, & \mbox{if} \,\, E(y)=4,\\
 \left(k_1^2D_1\right)^{m/2}+(-1)^{n/6+1}\left(3k_2^2\right)^{n/2}=1, & \mbox{if} \,\, E(y)=6,\\
\end{cases}
$$
where $D_1\neq1,\, D_2\neq1$ are square-free integers and $k_1,\, k_2$ are nonzero integers.  By Catalan's conjecture, we know that only two of these equations have non-trivial solutions, and these two equations are $\left(k_1^2D_1\right)^{m/2}-\left(k_2^2D_2\right)^{n/2}=1$ and $\left(k_1^2D_1\right)^{m/2}-k_2^n=1$. Then all non-trivial solutions of these two equations are
$$k_1^2D_1=\pm3,\, k_2^2D_2=2,\, m/2=2,\, n/2=3$$
and
$$k_1^2D_1=\pm3,\, k_2=2,\, m/2=2,\, n=3,$$
respectively. So we obtain
\[
\left(x,\, y,\, m,\, n\right)=\left(\pm\sqrt{\pm3},\, \pm\sqrt{2},\, 4,\, 6\right),\,\left(\pm\sqrt{\pm3},\, 2\omega,\, 4,\, 3\right),\,\left(\pm\sqrt{\pm3},\, 2\overline{\omega},\, 4,\, 3\right).
\]

{\bf Case 2:} $E(x)=3$. By Proposition \ref{p4}, equation \eqref{e9} becomes
$$
\begin{cases}
 k_1^m-\left(k_2^2D\right)^{n/2}=1, &  \mbox{if}\, \, E(y)=2, \\
 k_1^m-k_2^n=1, & \mbox{if} \,\, E(y)=3,\\
 k_1^m+(-1)^{n/4+1}\left(2k_2^2\right)^{n/2}=1, & \mbox{if} \,\, E(y)=4,\\
 k_1^m+(-1)^{n/6+1}\left(3k_2^2\right)^{n/2}=1, & \mbox{if} \,\, E(y)=6,\\
\end{cases}
$$
where $D\neq1$ is a square-free integer and $k_1,\, k_2$ are nonzero integers. By Catalan's conjecture, we know that these equations have no non-trivial solutions.

{\bf Case 3:} $E(x)=4$. By Proposition \ref{p4}, equation \eqref{e9} becomes
$$
\begin{cases}
(-1)^{m/4}\left(2k_1^2\right)^{m/2}-\left(k_2^2D\right)^{n/2}=1, &  \mbox{if}\, \, E(y)=2, \\
(-1)^{m/4}\left(2k_1^2\right)^{m/2}-k_2^n=1, & \mbox{if} \,\, E(y)=3,\\
(-1)^{m/4}\left(2k_1^2\right)^{m/2}+(-1)^{n/4+1}\left(2k_2^2\right)^{n/2}=1, & \mbox{if} \,\, E(y)=4,\\
(-1)^{m/4}\left(2k_1^2\right)^{m/2}+(-1)^{n/6+1}\left(3k_2^2\right)^{n/2}=1, & \mbox{if} \,\, E(y)=6,\\
\end{cases}
$$
where $D\neq1$ is a square-free integer and $k_1,\, k_2$ are nonzero integers. By Catalan's conjecture, we know that these equations have no non-trivial solutions.

{\bf Case 4:} $E(x)=6$. By Proposition \ref{p4}, equation \eqref{e9} becomes
$$
\begin{cases}
 (-1)^{m/6}\left(3k_1^2\right)^{m/2}-\left(k_2^2D\right)^{n/2}=1, &  \mbox{if}\, \, E(y)=2, \\
 (-1)^{m/6}\left(3k_1^2\right)^{m/2}-k_2^n=1, & \mbox{if} \,\, E(y)=3,\\
 (-1)^{m/6}\left(3k_1^2\right)^{m/2}+(-1)^{n/4+1}\left(2k_2^2\right)^{n/2}=1, & \mbox{if} \,\, E(y)=4,\\
 (-1)^{m/6}\left(3k_1^2\right)^{m/2}+(-1)^{n/6+1}\left(3k_2^2\right)^{n/2}=1, & \mbox{if} \,\, E(y)=6,\\
\end{cases}
$$
where $D\neq1$ is a square-free integer and $k_1,\, k_2$ are nonzero integers. By Catalan's conjecture, we know that these equations have no non-trivial solutions.
\end{prof}

\begin{theorem}\label{th6}
If the eigenvalues of $X$ and $Y$ are quadratic algebraic integers, then the following statements hold.
\begin{enumerate}
\item[\rm1)]\, If $X^m$ is a scalar matrix, then all non-trivial solutions of equation \eqref{e7} in $M_2(\mathbb{Z})$ are given by
$$
\left\{\left(X,\, Y,\, 4,\, 6\right): X,\, Y\in M_2(\mathbb{Z}),\, tr(X)=tr(Y)=0,\, \det(X)=\pm3,\, \det(Y)=-2\right\}
$$
and
$$
\left\{\left(X,\, Y,\, 4,\, 3\right): X,\, Y\in M_2(\mathbb{Z}),\, tr(X)=0,\, tr(Y)=-2,\, \det(X)=\pm3,\, \det(Y)=4\right\};
$$
\item[\rm2)]\, If $X^m$ is not a scalar matrix, then $XY=YX$.
\end{enumerate}
\end{theorem}
\begin{prof}
Assume that $x_s,\,y_s,\,s=1,\,2$ are the eigenvalues of $X$ and $Y$, respectively.

1)\, In this case, $x_1^m=x_2^m$ is an integer. From Lemma \ref{l7}, it follows that $(x_s,\, y_s,\, m,\, n),\,s=1,\,2$ are non-trivial solutions of equation \eqref{e9}. By Lemma \ref{l9}, we have
$$\left(x_s,\, y_s,\, m,\, n\right)\in\left\{\left(\pm\sqrt{\pm3},\, \pm\sqrt{2},\, 4,\, 6\right),\, \left(\pm\sqrt{\pm3},\, 2\omega,\, 4,\, 3\right),\, \left(\pm\sqrt{\pm3},\, 2\overline{\omega},\, 4,\, 3\right)\right\},\, s=1,\,2,$$
where $\omega=(-1+\sqrt{-3})/2$. By Lemma \ref{le5} and a direct computation, we can get all non-trivial solutions of equation \eqref{e7} in $M_2(\mathbb{Z})$ in this case, which are given in the theorem.

2)\, In this case, we have $x_1^m\neq x_2^m$. Since $(x_s,\, y_s,\, m,\, n),\,s=1,\,2$ are non-trivial solutions of equation \eqref{e9}, it follows that $x_1,\, x_2,\,y_1,\,y_2$ are quadratic algebraic integers in the same quadratic field. Let $K$ denote this quadratic field and let $\mathcal{O}_K$ be its ring of integers. Let $X=\begin{pmatrix} a_1 & b_1 \\ c_1 & d_1 \end{pmatrix}$ and $Y=\begin{pmatrix} a_2 & b_2 \\ c_2 & d_2 \end{pmatrix}$. Since $x_1,\, x_2,\,y_1,\,y_2\in \mathcal{O}_K\backslash\mathbb{Z}$, we have $b_1c_1\neq0$ and $b_2c_2\neq0$. By Lemma \ref{le5} and the assumption $X^m-Y^n=I$, we get the following three identities.
\begin{subequations}
\begin{align}
&\left(a_1-d_1\right)\cdot\frac{x_1^m-x_2^m}{x_1-x_2}=\left(a_2-d_2\right)\cdot\frac{y_1^n-y_2^n}{y_1-y_2} \label{a}\\
&\frac{y_1^n-y_2^n}{y_1-y_2}=\frac{b_1}{b_2}\cdot\frac{x_1^m-x_2^m}{x_1-x_2}\label{b} \\
&\frac{y_1^n-y_2^n}{y_1-y_2}=\frac{c_1}{c_2}\cdot\frac{x_1^m-x_2^m}{x_1-x_2} \label{c}
\end{align}
\end{subequations}
By \eqref{a} and \eqref{b}, we have
\begin{equation}\label{ea}
\left(a_1-d_1\right)b_2=\left(a_2-d_2\right)b_1.
\end{equation}
By \eqref{a} and \eqref{c}, we get
\begin{equation}\label{eb}
\left(a_1-d_1\right)c_2=\left(a_2-d_2\right)c_1.
\end{equation}
By \eqref{b} and \eqref{c}, we obtain
\begin{equation}\label{ec}
b_1c_2=b_2c_1.
\end{equation}
From \eqref{ea}, \eqref{eb} and \eqref{ec}, we conclude that $XY=YX$.
\end{prof}

Theorem \ref{th6} tells us that it is sufficient to study the solvability of the Catalan's matrix equation \eqref{e7} in $C(A)$, where $A=\begin{pmatrix} a & b \\ c & 0 \end{pmatrix}\in M_2\left(\mathbb{Z}\right)$ is a given matrix such that $bc\neq0$ and $\gcd(a,\, b,\, c)=1$. By Theorem \ref{th4}, we have the following corollary.
\begin{corollary}\label{c4}
Let $A=\begin{pmatrix} a & b \\ c & 0 \end{pmatrix}\in M_2\left(\mathbb{Z}\right)$ be a given matrix such that $bc\neq0$ and $\gcd(a,\, b,\, c)=1$. Let $K=\mathbb{Q}\left(\sqrt{a^2+4bc}\right)$ and let $\mathcal{O}_K$ be its ring of integers. Then the following statements hold.
\begin{enumerate}
\item[\rm1)]\,If $a^2+4bc$ is a square, then equation \eqref{e7} has no non-trivial solutions in $C(A)$;
\item[\rm2)]\,If $a^2+4bc$ is not a square and $D$ is the unique square-free integer such that $a^2+4bc=k^2D$ for some $k\in\mathbb{N}$, then equation \eqref{e7} has a non-trivial solution in $C(A)$ if and only if equation \eqref{e9} has a non-trivial solution $\left(x,\, y,\, m,\, n\right)$ in $\mathcal{O}_K$ such that $x,\, y$ can be written in the form
$$
\frac{s+t\sqrt{D}}{2},\quad s,\, t\in\mathbb{Z},\, k\mid t.
$$
\end{enumerate}
\end{corollary}

Hence, from Theorems \ref{th5}, \ref{th6} and Corollary \ref{c4}, we conclude that the solvability of the Catalan's matrix equation \eqref{e7} in $M_2\left(\mathbb{Z}\right)$ can be reduced to the solvability of the Catalan's matrix equation \eqref{e7} in $C(A)$, and finally to the solvability of the Catalan's equation \eqref{e9} in quadratic fields. However, the  solvability of the Catalan's equation in quadratic fields is unsolved. We leave this as an open question.

Let $K$ be a quadratic field and $\mathcal{O}_K$ its ring of integers. Let $D$ be the unique square-free integer such that $K=\mathbb{Q}\left(\sqrt{D}\right)$. From a given non-trivial solution of equation \eqref{e9} in $\mathcal{O}_K$, we can construct some classes of $2\times 2$ matrices such that equation \eqref{e7} has non-trivial solutions in these classes. Assume that
$$
\left(\frac{s_1+t_1\sqrt{D}}{2}\right)^m-\left(\frac{s_2+t_2\sqrt{D}}{2}\right)^n=1
$$
 is a given non-trivial solution of equation \eqref{e9} in $\mathcal{O}_K$, where $s_i,\, t_i\in\mathbb{Z},\, i=1,\, 2$. Let $a,\, b,\, c$ be integers and $k$ a positive integer such that $a^2+4bc=k^2D$, $k\mid t_1$, $k\mid t_2$, $bc\neq0$ and $\gcd\left(a,\, b,\, c\right)=1$. Indeed, such $a,\, b,\, c,\, k$ exist. If $D\equiv1\pmod4$ and $D=1+4t$ for some $t\in\mathbb{Z}$, then $\left(a,\, b,\, c,\,k\right)=\left(1,\, t,\, 1,\, 1 \right)$ satisfies the above conditions. If $D\equiv2,\, 3\pmod4$, then $\left(a,\, b,\, c,\, k\right)=\left(0,\, D,\, 1,\, 2\right)$ satisfies the above conditions. From the proof of Theorem \ref{th4}, it follows that
$$
\begin{pmatrix} \frac{s_1+\frac{t_1}{k}a}{2} & \frac{t_1}{k}b \\ \frac{t_1}{k}c & \frac{s_1-\frac{t_1}{k}a}{2} \end{pmatrix}^m-\begin{pmatrix} \frac{s_2+\frac{t_2}{k}a}{2} & \frac{t_2}{k}b \\ \frac{t_2}{k}c &  \frac{s_2-\frac{t_2}{k}a}{2} \end{pmatrix}^n=I
$$
are non-trivial solutions of equation \eqref{e7}, and the corresponding matrix classes are $C(A)$, where $A=\begin{pmatrix} a & b \\ c & 0\end{pmatrix}$. Next, we give an example to illustrate how to construct non-trivial solutions of the Catalan's matrix equation in $M_2\left(\mathbb{Z}\right)$ from a given equality in this manner.
\begin{example}\rm
Let $m,\, n\geq3$ be integers such that $m\equiv1\pmod6$ and $n\equiv-1\pmod6$. Then we have
$$
\left(\frac{1-\sqrt{-3}}{2}\right)^m-\left(\frac{-1+\sqrt{-3}}{2}\right)^n=1.
$$
From this equality, we can construct some classes of $2\times 2$ matrices such that equation \eqref{e7} has non-trivial solutions in these classes. Let $a,\, b,\, c$ be integers such that $a^2+4bc=-3$, $bc\neq0$ and $\gcd\left(a,\, b,\, c\right)=1$. Then
$$
\begin{pmatrix} \frac{1-a}{2} & -b \\ -c & \frac{1+a}{2} \end{pmatrix}^m-\begin{pmatrix} \frac{-1+a}{2} & b \\ c &  \frac{-1-a}{2} \end{pmatrix}^n=I
$$
are non-trivial solutions of equation \eqref{e7}, and the corresponding matrix classes are $C(A)$, where $A=\begin{pmatrix} a & b \\ c & 0\end{pmatrix}$.
\end{example}

\end{document}